\ifx\documentclass\undefined
\documentstyle[12pt]{article}
\else
\documentclass[12pt]{article}
\fi

\usepackage{color}

\usepackage{ascmac}

\newcommand{\beq}{\begin{eqnarray*}}
\newcommand{\eeq}{\end{eqnarray*}}
\makeatletter
\renewcommand{\theequation}{\thesection.\arabic{equation}}
\@addtoreset{equation}{section}
\def\eqnarray{%
\stepcounter{equation}%
\let\@currentlabel=\theequation
\global\@eqnswtrue
\global\@eqcnt\z@
\tabskip\@centering
\let\\=\@eqncr
$$\halign to \displaywidth\bgroup\@eqnsel\hskip\@centering
$\displaystyle\tabskip\z@{##}$&\global\@eqcnt\@ne
\hfil$\displaystyle{{}##{}}$\hfil
&\global\@eqcnt\tw@$\displaystyle\tabskip\z@{##}$\hfil
\tabskip\@centering&\llap{##}\tabskip\z@\cr}
\makeatother
\newtheorem{theorem}{Theorem}[section]

\newtheorem{remark}{Remark}[section]

\newsavebox{\toy}
\savebox{\toy}{\framebox[0.65em]{\rule{0cm}{1ex}}}
\newcommand{\QED}{\usebox{\toy}}
\def\nlni{\par\ifvmode\removelastskip\fi\vskip\baselineskip\noindent}

\begin{document}
%%%%%%% DOUBLE SPACED %%%%%%%%
\setlength{\baselineskip}{15pt}
\title{
On the 
asymptotic behavior of Bessel-like diffusions
}
\author{
Yuuki Shimizu
\thanks{
Department of Mathematics,
Gakushuin University,
1-5-1, Mejiro, Toshima-ku, Tokyo, 171-8588, Japan.
e-mail :
spkx3cw9@snow.ocn.ne.jp 
}
\and
Fumihiko Nakano
\thanks{
Department of Mathematics,
Gakushuin University,
1-5-1, Mejiro, Toshima-ku, Tokyo, 171-8588, Japan.
e-mail : 
fumihiko@math.gakushuin.ac.jp
}
}
%\date{最終更新日：}
\maketitle
%第一ページの番号を消す
%\thispagestyle{empty}
%%%%%%% ABSTRACT %%%%%%%%%%%%%
\begin{abstract}
We 
derive the asymptotic behavior of the transition probability density of the Bessel-like diffusions for 
``dimension"   
$\rho = 0$. 
\end{abstract}

%Mathematics Subject Classification (2000): 82B44, 81Q10

%\tableofcontents
%%%%% INTRODUCTION %%%%%%%%%%%%%%%%%%%%%%%%%%%%%%%%%
\section{Introduction}
\subsection{Background}
Let 
$\rho > 0$.
A Bessel process 
of dimension 
$\rho$ 
is a diffusion process on 
$[0, \infty)$ 
with generator 
\beq
{\cal L}_{\rho}
:=
\frac 12
\left(
\frac {d^2}{d x^2}
+
\frac {\rho-1}{x}
\frac {d}{dx}
\right), 
\quad
x > 0.
\eeq
If 
the origin is a regular boundary
(i.e., 
$0 < \rho < 2$), we impose the reflecting boundary condition. 
Then 
the transition probability density with respect to the speed measure
$m_{\rho}(dx) = 
2 x^{\rho-1} dx$
is
\beq
P_{\rho} (t,x,y)
:=
\frac {1}{2t}
(xy)^{- \nu}
\exp 
\left(
- \frac {x^2+y^2}{2t}
\right)
I_{\nu}
\left(
\frac {xy}{t}
\right)
\eeq
where 
$I_{\nu}$ 
is the modified Bessel function and 
$\nu := \frac {\rho}{2}-1$.
We thus have 
\beq
P_{\rho}(t; x, y)
\sim
\frac {1}
{
2^{\rho/2} 
\Gamma \left( 
\rho/2
%\frac {\rho}{2}
\right)
}
\cdot
\frac {1}
{
t^{\rho/2}
}, 
\quad
t \to \infty.
\eeq
Here and henceforth 
we denote by 
$f \sim g$ 
if 
$\lim \frac fg = 1$.
In this paper
we consider a diffusion process on 
$[0, \infty)$ 
with generator : 
\beq
{\cal L}
&:=&
\frac 12
\left(
\frac {d^2}{dx^2}
+
b(x) 
\frac {d}{dx}
\right), 
\quad
x > 0, 
\eeq
where 
$b \in L^1_{loc}[0, \infty)$ 
so that the left boundary 
$0$ 
is regular where the reflecting boundary condition is imposed. 
We  
assume that 
${\cal L}$
is asymptotically equal to the generator of the Bessel process : \\

\noindent
{\bf Assumption }
\beq
b(x)
& = &
\frac {\rho-1 + \epsilon(x)}{x} + \eta(x), 
\quad
x \ge 1
\eeq
{\it 
where 
$\rho \in {\bf R}$, 
$\lim_{x \to \infty} \epsilon (x) = 0$, 
$\eta \in L^1_{loc} [0, \infty)$ 
such that the following limit exists. 
$A:=\lim_{x \to \infty} \int_1^x \eta(u) du \in {\bf R}$. 
}
\\
\noindent
Assumption  
implies that the function  
\beq
W(x) := 
\exp \left(
\int_1^x b(u) du
\right), 
\quad
x > 0
\eeq
varies regularly at 
$\infty$ 
with index 
$\rho - 1$ ; 
that is 
\beq
\lim_{x \to \infty}
\frac {W(\lambda x)}{W(x)} = \lambda^{\rho-1}, 
\eeq
for any 
$\lambda > 0$.
We denote by 
$R_{\alpha}(\infty)$
(resp. $R_{\alpha}(0)$)
the totality of regularly varying functions at infinity
(resp. zero) with index
$\alpha$. 
Our aim is 
to study the asymptotic behavior of the transition probability  density of this diffusion as 
$t \to \infty$.
The answer is known for 
$\rho \ne 0$
\cite{KJM, OJM} 
which we recall in Subsection 1.2. 
%
%%%%%
\subsection{Known Results}
Set 
\beq
W(x)
&=&
\exp
\left(
\int_1^x b(u) du
\right),
\quad
x > 0
\\
s(x ; c) 
&:=&
\int_c^x \frac {du}{W(u)}, 
\quad
m(x)
=
2 \int_0^x W(u) du
\eeq
which leads to the canonical form 
${\cal L} = 
\frac {d}{dm(x)}
\frac {d}{ds(x)}$.
Let
$p(t; x, y)$ 
be the transition probability density with respect to 
$m(dx)$ 
which is equal to the Laplace transform of the spectral function 
$\sigma$. 
\beq
p(t ; 0, 0)
&=&
\int_{[0, \infty)}
e^{- \lambda t} 
d \sigma (\lambda), 
\quad
t > 0. 
\eeq
Let 
\beq
G_s (x,y)
=
\int_0^{\infty}
e^{- s t} 
p(t ; x, y) dt, 
\quad
s > 0
\eeq
be Green's function.
Then  
$h(s) := G_s (0,0)$
satisfies
\beq
h(s)
=
\int_{[0, \infty)}
\frac 
{d \sigma(\xi)}
{s + \xi}, 
\quad
s > 0, 
\eeq
and  
$h$ 
is the characteristic function associated to 
$\widetilde{m}(x) := m(s^{-1}(x))$
by Krein's correspondence
\cite{Ko-W}.
When 
$\rho \ne 0$, 
the answer to our question is : 
%

%
%%%%%%%%%%%%%%%
\begin{theorem}
\label{1}
\mbox{}\\
{\bf (1) (\cite{KJM} Theorem 4.2）}
If 
$\rho > 0$, 
\beq
p(t; x, y)
\sim
\frac {1}
{
2^{\rho/2} 
\Gamma \left(\rho/2 \right)
}
\cdot
\frac {1}
{
\sqrt{t}
W(\sqrt{t})
}, 
\quad
t \to \infty.
\eeq
{\bf (2) (\cite{OJM} Theorem 5.1）}
If 
$\rho < 0$,
\beq
p(t ; 0, 0)
- 
\frac {1}{m(\infty)}
&\sim&
\frac {1}{m(\infty)^2}
\frac {2^{\frac {\rho}{2}+1}}
{
| \rho |
\Gamma \left(
(2 - \rho)/2
\right)
}
\sqrt{t}W(\sqrt{t}), 
\quad
t \to \infty.
\eeq
\end{theorem}
We also 
recall the following result which is an important ingredient of the proof of our main theorem. 
Let 
$h^* (s)
=
(s h(s))^{-1}$
be the dual of 
$h$ 
which is the characteristic function associated to 
$\widetilde{m}^{-1}(x)$
\cite{Ko-W}.
Let 
$\sigma^*$ 
be the corresponding spectral function.
%
%
%%%%%%%%%%%%%%%%%%%%%%%%
\begin{theorem}
\label{2}
{\bf (\cite{OJM}, Proposition 5.1)}\\
If 
$\rho < 2$, 
\beq
\sigma^* (\lambda)
\sim
\frac {
2^{\frac {\rho}{2} + 1}
}
{
(2 - \rho)
\Gamma 
\left(
\frac {2 - \rho}{2}
\right)^2
}
\sqrt{\lambda}
W
\left(
\frac {1}{\sqrt{\lambda}}
\right), 
\quad
\lambda \to + 0.
\eeq
\end{theorem}
We note  
that 
Theorem \ref{2} is valid even for 
$\rho = 0$. 
%
%
%%%%%%%%%%%%%%%%
\subsection{Results in this paper}
In this paper, 
we consider the case 
$\rho = 0$. 
Then we could have both 
$m(+\infty) = \infty$ 
and 
$m(+\infty) < + \infty$. 
Let $m_{\infty} := m(+\infty)$.
Since 
$\sigma(+0) = 
1/m_{\infty}$, 
$m_{\infty}< \infty$ 
implies 
$\sigma(+0) > 0$ 
and 
$p(t ; 0, 0)
\stackrel{t \to \infty}{\to}
1/m_{\infty}$. 
\begin{theorem}
\label{main}
If 
$\rho = 0$ 
and 
$m_{\infty} = \infty$, 
\beq
p(t; x, y)
&\sim&
\frac {1}{m (\sqrt{t})}, 
\quad
t \to \infty.
\eeq
\end{theorem}
\begin{theorem}
\label{main2}
If 
$\rho = 0$ 
and 
$m_{\infty} < \infty$, 
\beq
p(t ; 0, 0) - \frac {1}{m_{\infty}}
& \sim &
\frac {1}{m_{\infty}^2}
\left(
m_{\infty} - 
m
\left(\sqrt{t}\right)
\right), 
\quad
t \to \infty.
\eeq
\end{theorem}
\begin{remark}
\label{unify}
To summarize 
the statements in 
\cite{KJM} Theorem 4.1, 
\cite{OJM} Theorem 5.1
and 
Theorems \ref{main}, \ref{main2}, 
we have 
\\
(1)
$\rho \ge 0$, 
$m(+\infty) = \infty$ : 
\begin{equation}
p(t; x, y)
\sim
\frac {1}
{
2^{\rho/2}
\Gamma\left(
\rho/2 + 1
\right)
}
\cdot
\frac {1}
{
m (\sqrt{t})
},
\quad
t \to \infty.
\label{infinite}
\end{equation}
(2)
$\rho \le 0$, 
$m(+\infty) < \infty$ : 
\begin{equation}
p(t ; 0, 0) - \frac {1}{m_{\infty}}
\sim 
%\frac {1}{m_{\infty}}
\frac {
1
%2^{\frac {\rho}{2}}
}
{
2^{|\rho|/2}
\Gamma \left( 
|\rho|/2 + 1
%\frac {2 - \rho}{2} \right)^2
\right)
}
%\Gamma \left( 1 - \frac {\rho}{2} \right)
\cdot
\frac {1}{m_{\infty}}
\left(
1 - 
\frac {
m
\left(
\sqrt{t}
\right)
}
{m_{\infty}
}
\right), 
\quad
t \to \infty.
\label{finite}
\end{equation}
\end{remark}
In Section 2, 
we prove Theorems \ref{main}, \ref{main2} and apply them to some concrete examples. 
A strategy 
of the proof is to study the behavior of the following quantities in the arranged order, using Theorem \ref{2} and Tauberian theorems. 
\beq
\sigma^* (\lambda)
\to
h^*(s)
\to
h(s) 
=
\frac {1}{s h^*(s)}
\to
\sigma(\lambda)
\eeq
In Section 3, 
we shall quote some Tauberian Theorems used frequently in this paper. 
%
%%%%%%%%%%%%%%%%%%%%%%%%
\section{Proof of Theorems}
%%%%%%%%%
\subsection{Proof of Theorem \ref{main}
}
First of all, 
by a property of the regularly varying functions \cite{BGT} we have 
\beq
m(x)
&=&
2 \int_0^x W(u) du
\sim
\frac {2}{\rho}\,
x W(x), 
\quad
x \to \infty.
\eeq
Applying it to 
Theorem \ref{1} 
yields 
(\ref{infinite})
in Remark \ref{unify} for 
$\rho > 0$. 
\\

\noindent
{\it Proof of Theorem \ref{main}}\\
By the argument in 
\cite{KJM} Corollary 5.3, 
\[
p(t, x, y) \sim p(t, 0, 0), 
\quad
t \to \infty
\]
so that we may suppose 
$x = y = 0$.
\cite{OJM} 
Proposition 5.1
($\rho = 0$) 
implies 
\beq
\sigma^* (\lambda)
\sim
\sqrt{\lambda}
W 
\left(
\frac {1}{\sqrt{\lambda}}
\right)
\in R_1(0), 
\quad
\lambda \downarrow 0.
\eeq
Thus 
\cite{OJM} 
Proposition 5.1
($\rho = 0$) and 
Theorem \ref{Tauberianh} 
($\alpha = 1$, $n=1$) 
below yield
\beq
(-1)\cdot
\frac {d}{ds}h^*(s)
\sim
s^{-2} \sigma^*(s)
\sim
s^{- \frac 32} 
W \left(
\frac {1}{\sqrt{s}}
\right),
\quad
s \to +0.
\eeq
On the other hand, 
by the definition of 
$m$, 
\beq
\frac {d}{ds}
m
\left(
\frac {1}{\sqrt{s}}
\right)
&=&
-s^{- \frac 32}
W \left(
\frac {1}{\sqrt{s}}
\right).
\eeq
Therefore
\begin{equation}
- \frac {d}{ds}
h^* (s)
\sim
- \frac {d}{ds}
m \left(
\frac {1}{\sqrt{s}}
\right), 
\quad
s \to +0.
\label{bibun}
%\quad\cdots (*)
\end{equation}
Since 
$m(+\infty) = \infty$, 
we may apply 
de l'Hospital's theorem to have 
\beq
h^* (s)
\sim
m \left(
\frac {1}{\sqrt{s}}
\right), 
\quad
s \to +0.
\eeq
Using 
$h^* (s) = (s h(s))^{-1}$, 
\begin{equation}
h(s) \sim
\frac {1}{
s 
m \left(
\frac {1}{\sqrt{s}}
\right)
}, 
\quad
s \to +0.
\label{h}
\end{equation}
Note that, 
by 
\cite{BGT} Proposition 1.5.9a 
and by the fact that 
$l(x) := x W(x)$ 
is slowly varying at infinity, 
$f(s) := m\left( \frac {1}{\sqrt{s}} \right)$ 
is slowly varying at 
$0$. 
By Theorem \ref{Tauberianh} 
($\alpha=0$, $n=0$) below, 
\beq
\sigma(s)
\sim
\frac {1}
{
m \left(
\frac {1}{\sqrt{s}}
\right)
}, 
\quad
s \to +0.
\eeq
Thus  
Theorem \ref{Tauberian} 
completes the proof. 
\QED
\\
\begin{remark}
There is 
another argument starting with 
(\ref{h}).
Using  
$h(s) = \int_0^{\infty} e^{-st} p(t ; 0, 0) dt$, 
Theorem \ref{Tauberian} 
implies 
\beq
\int^t p(s; 0, 0) ds
& \sim &
\frac {t}{m(\sqrt{t})}, 
\quad
t \to \infty.
\eeq
Since 
$p(t, 0,0)$ 
is monotone as a function of 
$t$, 
monotone density theorem 
(\cite{BGT}, Theorem 1.7.2)
yields 
\beq
p(t; 0, 0) & \sim & \frac {1}{m(\sqrt{t})}, 
\quad
t \to \infty.
\eeq
\end{remark}
%
%%%%%%%%%%%%%%%%%%%%%%%%%%%%%%%%%%%%
\subsection{Proof of Theorem \ref{main2}}
%%%%%%
%
We first 
derive 
(\ref{finite}) 
in Remark \ref{unify} for 
$\rho < 0$. 
Set 
\beq
m(x)  &=& 
2 \int_0^x W(u) du, 
\quad
m_{\infty} = 2 \int_0^{\infty} W(u) du, 
\\
s(x) &=& 
\int_0^x \frac {1}{W(u)} du.
\eeq
Let 
${\cal L}^{\bullet}$
be the dual operator of 
${\cal L}$ 
\beq
{\cal L}^{\bullet}
&:=&
\frac 12
\left(
\frac {d^2}{dx^2}
-
b(x) 
\frac {d}{dx}
\right)
\eeq
and let 
$m^{\bullet}$, $s^{\bullet}$ 
be the corresponding speed measure and the scale function, respectively.
Then 
\beq
m^{\bullet}(x) &=&
2 \int_0^x \frac {1}{W(u)} du = 2 s(x)
\\
s^{\bullet}(x) &=&
\int_0^x W(u) du 
=
\frac 12 m(x)
\eeq
so that 
$l^{\bullet} 
:= h^{\bullet}(+0)
= s^{\bullet}(+\infty) 
= \frac 12 m_{\infty}$.
Since 
$h^*(s) = 2 h^{\bullet}(s)$
\cite{OJM}, 
we have 
\begin{equation}
l^* 
:= h^*(+0) 
= 2 h^{\bullet}(+0) 
= 2 l^{\bullet}
=m_{\infty}.
\label{equation}
\end{equation}
Thus
\beq
h(s)
&=&
\frac {1}{s h^*(s)}
\sim
\frac {1}{s m_{\infty}},
\quad
s \to +0.
\eeq
and by 
Theorem \ref{Tauberianh}
($\alpha=0$, $n=0$, $A=m_{\infty}^{-1}$) below, 
%we have 
%
\beq
\sigma(\lambda)
 \sim 
\frac {1}{m_{\infty}}, 
\quad
\lambda \to + 0. 
\eeq
Since 
\beq
\int_x^{\infty} W(u) du
& \sim &
\frac {1}{| \rho |} x W(x), 
\quad
x \to \infty
\eeq
by 
\cite{BGT} Proposition 1.5.10, 
we have 
\beq
m_{\infty} - m(x)
& \sim &
\frac {2}{| \rho |} x W(x)
\in 
R_{\rho}(0)
\eeq
which, together with 
\cite{OJM} Theorem 5.1, 
yields 
(\ref{finite}) 
in 
Remark \ref{unify}. 
\\

\noindent
{\it Proof of Theorem \ref{main2}}
We note that, by 
\cite{BGT} Proposition 1.5.9b, 
$g(s) := 
m_{\infty} - m \left( \frac {1}{\sqrt{s}} \right)$
is slowly varying at $0$. 
Owing to 
Theorem \ref{Tauberianh}, 
it suffices to show the following equation. 
\beq
h(s) - \frac {1}{s m_{\infty}}
& \sim &
\frac {1}{s m_{\infty}^2}
\left(
m_{\infty} - m \left( \frac {1}{\sqrt{s}} \right)
\right), 
\quad
s \to +0
\eeq
which is equivalent to 
\beq
\frac {1}{h^* (s)} - \frac {1}{m_{\infty}}
& \sim &
\frac {1}{m_{\infty}^2} 
\left(
m_{\infty} - m \left( \frac {1}{\sqrt{s}} \right)
\right).
\eeq
By 
de l'Hospital's theorem, 
\beq
\frac {
\frac {1}{h^*(s)} - \frac {1}{m_{\infty}}
}
{
\frac {1}{m_{\infty}^2}
\left(
m_{\infty} - m \left( \frac {1}{\sqrt{s}} \right)
\right)
}
& \sim &
\frac {
\left( \frac {1}{h^* (s)} \right)'
}
{
\frac {1}{m_{\infty}^2}
\left(
- \frac {d}{ds}
m \left( \frac {1}{\sqrt{s}} \right)
\right)
}
\\
&=&
\frac {
-\frac {(h^{*})'(s)}{h^{*}(s)^2}
}
{
\frac {1}{m_{\infty}^2}
\left(
- \frac {d}{ds}
m \left( \frac {1}{\sqrt{s}} \right)
\right)
}
\\
& \sim &
\frac {
\frac {1}{m_{\infty}^2}
\frac {d}{ds}
m \left( \frac {1}{\sqrt{s}} \right)
}
{
\frac {1}{m_{\infty}^2}
\frac {d}{ds}
m \left( \frac {1}{\sqrt{s}} \right)
}
= 1, 
\quad
s \to +0
\eeq
which finishes the proof, where we used 
(\ref{bibun}), (\ref{equation}) 
in the last line. 
\QED
%
%%%%%%%%%%
\subsection{Example}
We applly 
Theorems \ref{main}, \ref{main2} 
to some examples. 
In what follows, 
$\eta \in L^1_{loc}[0, \infty)$ 
such that the limit  
$A:= 
\lim_{x \to \infty} \int_1^x \eta(u) du$
exists. 
\\
\noindent
{\bf Example 1}
\beq
b(x)
&=&
- \frac 1x 1(x \ge 1)
+
\eta (x), 
\eeq
Then we have 
\beq
p(t ;  x, y)
\sim
\frac {e^{-A}}{2}
\left(
\log \sqrt{t}
\right)^{-1}, 
\quad
t \to \infty.
\eeq
\noindent
{\bf Example 2}
\beq
b(x) 
&=&
\left(
- \frac 1x + \frac {\alpha}{x (\log x)^{\beta}}
\right)
1(x > 1)
+
\eta(x), 
\quad
\alpha \ne 0,
\quad
0 < \beta < 1.
\eeq
Note that the case 
$\beta > 1$
is reduced to Example 1.
Then 
\beq
p(t ; x, y)
& \sim &
\frac {\alpha}{2}
e^{-A}
(\log \sqrt{t})^{-\beta}
e^{
-\frac {\alpha}{1 -\beta}
(\log \sqrt{t})^{1 - \beta}
}, 
\quad
t \to \infty.
\eeq
%
%%%%%%
{\bf Example 3}
\beq
b(x)
&=&
\left(
- \frac 1x + 
\frac {\alpha}{x \log x}
\right)
1(x > e)
+
\eta(x)
\eeq
Then
\beq
&(1)&\quad \alpha > -1, 
\quad
p(t; x, y)
\sim
\frac {\alpha+1}{2} e^{-(A+1)}
(\log \sqrt{t})^{- (\alpha+1)}, 
\quad
t \to \infty
\\
&(2)&\quad \alpha = -1, 
\quad
p(t; x, y)
\sim
\frac 
{
e^{-(A+1)}
}{2} 
(\log \log \sqrt{t})^{-1}, 
\quad
t \to \infty
\\
&(3)&\quad \alpha < -1, 
\quad
p(t ; x, y) - \frac {1}{m_{\infty}}
\sim 
\frac {1}{m_{\infty}^2}
\frac {(-2)}{\alpha+1}
e^{A+1}
(\log \sqrt{t})^{\alpha +1}, 
\quad
t \to \infty.
\eeq
where 
$m_{\infty}
:=
2 \int_0^{\infty} 
\exp 
\left(
\int_1^u b(v) dv
\right)
du$. 
\\
\noindent
{\bf Example 4}\\
In general, 
given a function 
$m : [0, \infty) \to (0, \infty)$, 
such that 
$\lim_{t \to \infty}
\frac {m''(t)}{m'(t)} t = -1$, 
we can construct a corresponding generator 
${\cal L}$ 
such that 
$p(t ; x, y) \sim (m(\sqrt{t}))^{-1}$, 
$t \to \infty$. 
In fact, 
we can take 
\beq
b(x) :=
-\frac 1x
+
\frac {f''(\log x)}{f'(\log x)}
\cdot
\frac 1x
\eeq
where 
$f(x) := m(e^x)$.
%

%
%%%%%%%%%%%%%%%%%%%%%%%%%%%%%
\section{Appendix}
We recall 
some important facts from the theory of regularly varying functions 
\cite{BGT}, \cite{OJM}.
For a function 
$\sigma : [0, \infty) \to {\bf R}$
being of locally bounded variation and right-continuous, let 
\beq
\hat{\sigma}(\lambda)
&=&
\int_{[0, \infty)} e^{- \lambda x}d \sigma(x)
\\
H_n (\sigma, \lambda)
&:=&
\int_{[0, \infty)}
\frac {d \sigma(\xi)}
{(\lambda + \xi)^{n+1}}, 
\quad
n \ge 0
\eeq
be its Laplace transform, and the generalized Stieltjes transform, respectively. 
\begin{theorem}
\label{Tauberian}
Let 
$\rho \ge 0$
and 
$f \in R_{\alpha}(0)$. 
Then 
\beq
\sigma(x)
& \sim &
c 
f(x), 
\;
x \to \infty
\quad
\Longleftrightarrow 
\quad
\hat{\sigma}(\lambda)
\sim 
c 
\Gamma(\rho+1)
f \left( \frac {1}{\lambda} \right), 
\;
\lambda \to +0.
\eeq
\end{theorem}
\begin{theorem}
{\bf (Theorem 7.1 in \cite{OJM})}
\label{Tauberianh}
\\
Let 
$0 \le \alpha < n+1$, 
$A \ge 0$, 
and 
$\varphi \in R_{\alpha}(0)$.
Then
\beq
\sigma(\xi)
&\sim&
A \varphi (\xi), 
\;
\xi \to 0
\quad
\Longleftrightarrow
\quad
H_n(\sigma ; \lambda)
\sim 
A C_{n, \alpha} 
\varphi(\lambda) \lambda^{-n-1}, 
\;
\lambda \to 0
\eeq
where 
\[
C_{n, \alpha}
:=
\frac {
\Gamma(n+1 - \alpha) \Gamma(\alpha+1)
}
{
\Gamma(n+1)
}.
\]
\end{theorem}

\vspace*{1em}
\noindent {\bf Acknowledgement }
The authors 
would like to thank Professor Yuji Kasahara for his important and essential suggestions, kind advice and warm encouragement. 
The authors 
would also like to thank Professor Shin'ichi Kotani for pointing out Example 4, and the referee for many valuable  remarks. 
The second 
author(F.N.) is partially supported by 
JSPS KAKENHI Grant Number .26400145.

%This work is partially supported by 
%Grant-in-Aid for Scientific Research (C) no.26400145.

%%%%% REFERENCES %%%%%%%%%%%%%%%%%%%%%
%

\end{document}